\documentclass[12pt,a4paper]{article}
\usepackage[cp1251]{inputenc}
\usepackage{graphicx}
\usepackage{amssymb}
\usepackage{amsmath}
\usepackage{amsthm}
\usepackage{euscript}
\numberwithin{equation}{section}
\newtheorem{lemma}{Lemma}
\newtheorem{conjecture}{Conjecture}
\newtheorem{rem}{Remark}

\newtheorem{ex}{Example}

\newtheorem{corollary}{Corollary}

\newcommand{\N}[0]{\mathbb{N}}
\newcommand{\Z}[0]{\mathbb{Z}}

\newcommand{\be}[0]{\begin{equation}}
\newcommand{\ee}[0]{\end{equation}}
\newcommand{\bez}[0]{\begin{equation*}}
\newcommand{\eez}[0]{\end{equation*}}
\newcommand{\bl}[0]{\begin{lemma}}
\newcommand{\el}[0]{\end{lemma}}

\newcommand{\ep}[0]{$\hspace{\fill} \square$}
\newcommand{\paragraf}[1]{\par
\bigskip{\centerline{\bf #1}}\medskip}

\newcommand{\abs}[1]{\begin{quotation} {\small
 \centerline{{\bf Abstract}}\smallskip
#1}
\end{quotation}}
\author{I.~Shnurnikov\footnote{NRU HSE}}
\title{What is the number of decompositions of torus into given number of regions by unions of geodesics?}
\date{}
\begin{document}
\maketitle
\abs{ We prove some preliminary results concerning two questions of O.Karpenkov:

(1) What is the number of decompositions (up to $SL(2,\Z)$) of two--dimensional torus into given number $f$ of regions by unions of $n$ geodesics?

(2) On the plane there are $n$ circles not in general position, every pair of cicles has at least one common point. What is the set of all possible numbers of regions? 

}

\paragraf{Introduction}

Let us consider a flat two--dimensional torus $T$ (quotient space of real euclidean plane for an action of lattice --- abelian group with two generators). In a fixed homology basis on torus $T$ a closed oriented geodesic is defined (up to parallel shifting) by a pair of coprime integers. The matrix of changing from a homology basis to any other is an integer $2\times2$ matrix with determinant $\pm 1$. We shall consider arrangements of nonoriented geodesics up to changing homology basis. 
Let $f$ be the number of connected components of the complement in two--dimensional torus $T$ to the union of $n$ geodesics. 
The set $F(n)$ of all possible numbers $f$ for given $n>1$ is the following (see \cite{Shn11})
$$
F(n) = \{ n-1, n \} \cup \{ m \in \N \ | \ m \geq 2n-4 \}.
$$

Let $t_i$ be the number of intersection points, which are incident to $i$ geodesics of the arrangement. If not all geodesics of the arrangement are parallel, then $f=\sum_i (i-1) t_i$.
For example, geodesics of types $(1,0), (0,1), (k,1)$ form $k$ or $k+1$ regions if they intersect in one point or not.
If in arrangement of geodesics $\gamma_1, \dots, \gamma_n$ are at least two non-parallel, then
$$
f \leq \sum_{i<j} |\gamma_i \cap \gamma_j|
$$
where $|\gamma_i \cap \gamma_j|$ is the number of intersection points of non-parallel geodesics. For arrangements of general position the inequality turns to equality.

\paragraf{Question for torus}
In connection with the theory of high--dimensional chain fractions  O.Karpenkov asked:

What is the number of decompositions of two--dimensional torus into given number $f$ of regions by unions of $n$ geodesics?


\begin{lemma}
If two geodesics intersect in $k$ points in the two--dimensional flat torus, then we may change bases so that geodesics will be of type $(1,0)$ and $(x,k)$, where integer $1 \leq x \leq k-1$ is such that $gcd(x,k)=1$ and is defined uniquely up to change $x \leftrightarrow k-x$.
\end{lemma}

\begin{lemma}
For $f=n$ and $f=n-1$ there is a unique arrangement of $n$ geodesics in the two--dimensional torus which divides torus into $f$ regions.
The number $f=2n-4$ is realised as the number of regions by $n-1$ arrangements for $n \geq 7$ ( and at least by 8 arrangements for $n=6$). 
\end{lemma}
\proof 
Let us take $n-2$ geodesics of type $(1,0)$ a geodesic of type $(a,1)$ and a geodesic of type $(0,1)$, where $0 \leq a \leq n-2$, and all intersection points of the last two geodesics are incident to some of the first $n-2$ geodesics.
\ep

Let $m$ be the maximal number of parallel (homologically equal) geodesics in an arrangements.

\begin{lemma}
If $f \leq c n^{\tfrac 65}$, then $m \geq n- \frac fn +O(1)$,
for suitable positive coonstant $c$.

\end{lemma}

\begin{corollary}
For $f \leq c n^{\tfrac 65}$ in arrangement of geodesics almost all geodesics are homologically equal and so the number of arrangements which realize $f$ as the number of regions may be counted explicitly.
\end{corollary}

 \paragraf{Question for circles in the plane}

O.Karpenkov asked the following ``On the plane there are $n$ circles not in general position, every pair of cicles has at least one common point. What is the set of all possible numbers of regions for given $n$?'' 

For comlete solution of this problem one need to determine the possible number of tangent points of $a$ circles and special arrangements of $n-a$ lines for $n>ca^2$. 

Let us denote by $C_n$ the set of numbers $f$ which are formed  by $n$ circles in the plane not in general position such that every two circles have at least one common point. Let us denote by $L_n$ the set of numbers of regions in the plane, formed by $n$ distinct lines not in general position (without any requirements on intersection points). Let $m$ be the maximal number of circles, incident to one point.

\begin{lemma}
We have $C_n\supseteq L_n$.
\end{lemma}
\proof 
 Let us take any arrangement of lines in the plane with $f \in L_n$, make an inversion and get the suitable arrangement of $n$ circles with $f$ regions. 
\ep 

If $m=n$, i.e. all circles have one common point, then the number of regions $f \in L_n$.
\begin{lemma}
If $m=n-1$, then $f$ may be any number of the sets $L_{n-1}+2n-2$, $L_{n-1}+2n-3$, $L_{n-1}+2n-4$, and this list of possibilities is uncomplete (here we sum a number to every element of $L_{n-1}$). 

The numbers $3n-4, 4n-4 \in C_n$ and $3n-4, 4n-4 \notin L_{n}$,
We have
$$
L_n=\{n+1,2n,3n-3,3n-2,4n-8,4n-7,4n-6,4n-5,5n-15, \dots \}
$$ 

\end{lemma} 

\begin{conjecture}
The set $C_n$ contains all integers  between $\frac{n(n-1)}2+1$ and $n(n-1)+2$, which are the maximal elements of $L_n$ and $C_n$ correspondingly. 
\end{conjecture}

\end{document}